\documentclass[11pt]{amsart}
\usepackage{mathrsfs,amsthm,amssymb,amsmath,url}
\usepackage[latin1]{inputenc}
\def\cal{\mathscr}

\theoremstyle{plain}
    \newtheorem{thm}{Theorem}
    
    \newtheorem{lem}[thm]{Lemma}
    \newtheorem{prop}[thm]{Proposition}
    \newtheorem{cor}[thm]{Corollary}

\theoremstyle{definition}
    \newtheorem{ex}[thm]{Example}
    \newtheorem{defn}[thm]{Definition}
    
    \newtheorem{prob}[thm]{Problem}

\theoremstyle{remark}
    \newtheorem{rem}[thm]{Remark}

\newcommand{\nin}{\notin}

\newcommand{\rest}{\upharpoonright}

\renewcommand{\P}{{\mathscr P}}

\newcommand{\cl}[1]{\langle #1 \rangle}
 \DeclareMathOperator{\supp}{supp}
\DeclareMathOperator{\Cl}{Cl}
\DeclareMathOperator{\id}{id}\DeclareMathOperator{\proj}{Proj}

\renewcommand{\O}{{\mathscr O}}
\newcommand{\On}{{\mathscr O}^{(n)}}

\newcommand{\Oo}{{\mathscr O}^{(1)}}

\DeclareMathOperator{\Pol}{Pol} \DeclareMathOperator{\pol}{Pol}

\newcommand{\inv}{^{-1}}
\newcommand{\U}{{\mathscr U}}
\newcommand{\C}{{\mathscr C}}

\newcommand{\F}{{\mathscr F}}

\newcommand{\A}{{\mathscr A}}

\newcommand{\D}{{\mathscr D}}

\newcommand{\G}{{\mathscr G}}
\newcommand{\M}{{\mathscr M}}
\newcommand{\N}{{\mathbb N}}

\renewcommand{\S}{{\mathscr S}}

\newcommand{\un}{^{(n)}}

\newcommand{\uo}{^{(1)}}

\newcommand{\ida}{\hat{I}}

\newcommand{\To}{\rightarrow}

\newcommand{\restr}{\rest}
\renewcommand{\subset}{\subseteq}
\renewcommand{\supset}{\supseteq}

\newcommand{\abs}[1]{|#1|}
\newcommand{\bool}[1]{[\![\, #1\, ]\!]}

 \newtheorem{obss}[thm]{Observations}
\newcommand{\weg}[1]{ }

\newcommand{\xitem}[1]{\item[$(#1)$]}
\newcommand{\BE}{\begin{enumerate}}
\newcommand{\EE}{\end{enumerate}}

\title{Clones from ideals}
\author[M. Beiglb\"{o}ck]{Mathias Beiglb\"{o}ck}
\address{Algebra\\TU Wien\\Wiedner Hauptstrasse 8-10/104\\1040 Wien, Austria}
\email{mathias.beiglboeck@tuwien.ac.at}\urladdr{http://www.dmg.tuwien.ac.at/beigl/}
\author[M. Goldstern]{Martin Goldstern}
\address{Algebra\\TU Wien\\Wiedner Hauptstrasse 8-10/104\\1040 Wien, Austria}
\email{goldstern@tuwien.ac.at}\urladdr{http://www.tuwien.ac.at/goldstern/}
\author[L. Heindorf]{Lutz Heindorf}
\address{Mathematisches Institut\\FU Berlin\\Arnimallee 3\\14195 Berlin, Germany}
\email{heindorf@math.fu-berlin.de}
\author[M. Pinsker]{Michael Pinsker}
\address{Laboratoire de Math\'{e}matiques Nicolas Oresme\\
CNRS UMR 6139\\ Universit\'{e} de Caen\\14032 Caen Cedex\\
France}
\email{marula@gmx.at}\urladdr{http://dmg.tuwien.ac.at/pinsker/}

\thanks{The first and second author are grateful to the Austrian science foundation (FWF) for support under grant  P17627; the fourth author is grateful to the
Austrian science foundation for support under grant  P17812.}

\begin{document}
\maketitle
\begin{abstract}
    On an infinite base set~$X$, every ideal of subsets of~$X$ can
    be associated with the clone of those operations on~$X$ which
    map small sets to
    small sets. We continue earlier investigations on the position
    of such
    clones in the clone lattice.
\end{abstract}

\setcounter{section}{-1}
\section{Introduction}

\subsection{Clones}
Let~$X$ be an infinite set and denote the set of all~$n$-ary
operations on~$X$ by $\On$. Then $\O:=\bigcup_{n\geq 1}\O\un$ is the
set of all finitary operations on~$X$. A subset $\C$ of $\O$ is
called a \emph{clone} iff it contains all projections, i.e.\ for all
$1\leq k\leq n$ the function~$\pi^n_k\in\On$ satisfying
$\pi^n_k(x_1,\ldots,x_n)=x_k$, and is closed under composition. The
set of all clones on~$X$, ordered by set-theoretical inclusion,
forms a complete algebraic lattice $\Cl(X)$. The size of this
lattice is $2^{2^{|X|}}$, and it is known to be too complicated to
be ever fully described: For example, it contains all algebraic
lattices with no more than $2^{|X|}$ compact elements as complete
sublattices \cite{Pin07AlgebraicSublattices}. So the approach to an
unveiling of the structure of $\Cl(X)$ is to concentrate on more
tractable parts of it, such as large clones (e.g.\ dual atoms of the
lattice), or clones with specific properties (e.g.\ natural
intervals, clones closed under conjugation, etc.). A survey of
clones on infinite sets is \cite{GP06survey}.

\subsection{Ideal clones} In \cite{Ros76}, and before in \cite{Gav65} for countably infinite~$X$,
it was shown that there exist as many dual atoms
(``\emph{precomplete clones}'') in $\Cl(X)$ as there are clones,
suggesting that it is impossible to describe all of them (as opposed
to the clone lattice on finite~$X$, where the dual atoms are finite
in number and explicitly known \cite{Ros70}). Much more recently, a
new and short proof of this fact was given in
\cite{GS02clonesOnRegular}. It was observed that given an ideal~$I$
of subsets of~$X$, one can associate with it a clone $\C_I$
consisting of those operations which map small sets (i.e., products
of sets in~$I$) to small sets.
The authors then showed that prime
ideals correspond to precomplete clones, and that moreover the
clones induced by distinct prime ideals differ, implying that there
exist as many precomplete clones as prime ideals on~$X$; the latter
are known to amount to $2^{2^{|X|}}$.

The study of clones that arise in this way from ideals was pursued
in \cite{CH01}, for countably infinite~$X$. The authors concentrated
on the question of which ideals induce precomplete clones, and
provided a criterion for precompleteness.

\subsection{Precompleteness criteria}\label{sec:precompleteness}

We will consider a countable base set~$X$;  all our ideals $I$ will
be proper and properly contain the ideal of finite subsets of $X$.
(See~\ref{our.ideals} for an explanation.) Under those assumptions
on~$I$, which will be valid throughout this paper, it was shown in
\cite{CH01} that $\C_I$ is precomplete if and only if for all $A\nin
I$ there exists an operation~$f\in\C_I$ such that $f[A^n]=X$.

\subsubsection{Arity}
It is a  drawback of this
criterion  that, given a set $A\nin I$, we do not know in advance
the arity of the required function  mapping~$A$ onto~$X$.  It would be of much help if one had to
check only whether, say, a binary operation of $\C_I$ can map~$A$
onto~$X$. Unfortunately this is not the case: we will see in this paper that for every
$n\geq 1$, there exists an ideal~$I_n$ on a countably infinite base
set~$X$ such that $\C_{I_n}$ is precomplete, and its precompleteness
can be verified using  $(n+1)$-ary operations,
but not with $n$-ary operations.

\subsubsection{Regular ideals}
A better criterion would be one where precompleteness of $\C_I$ can
be read off the ideal~$I$ directly, without the use of the functions
in $\C_I$. We will obtain such a criterion for certain ideals~$I$,
on a countably infinite base set~$X$: Define $\hat{I}$ to consist of
all subsets~$A$ of~$X$ such that every infinite subset~$B$ of~$A$
contains an infinite set in~$I$. We will prove that if $X
\not\in\ida$, then $\C_I$ is precomplete iff $I=\ida$; so the only
ideals left to consider are those which satisfy $X\in \ida$, i.e.\
those for which $\ida$ equals the power set of $X$. Utilizing this
new criterion, we will once again construct $2^{2^{\aleph_0}}$
precomplete clones on a countably infinite base set, but for the
first time in ZF, i.e.\ explicitly and without the use of the Axiom
of Choice.

\subsubsection{A unary criterion}\label{ukrit}
In \cite[Problem D]{GP06survey}, it was asked whether there is
 a precompleteness criterion for $\C_I$ using unary operations
only.  The authors justified this question by
observing that every clone $\C_I$ is actually determined by its
unary operations: For a transformation monoid $\G\subseteq\Oo$, define
$\pol(\G)\subseteq\O$ to consist of those operations~$f$ satisfying
$f(g_1,\ldots,g_n)\in\G$ for all $g_1,\ldots,g_n\in\G$. Then we have
$\C_I=\pol(\C_I\uo)$. We will provide a criterion using only unary
functions in this paper.

\subsubsection{Extension to precomplete ideal clones} In \cite[Problem 1]{CH01}, the authors asked
whether every ideal clone could be extended to a precomplete ideal
clone. We will provide a positive answer to this question, and
moreover show that every precomplete clone extending an ideal clone
is itself an ideal clone. It follows that an ideal clone is
precomplete iff it is not properly contained in a precomplete ideal
clone.

\subsection{The mutual position of ideal clones}

Our investigations will not exclusively treat the precompleteness of
ideal clones, but also their mutual position.

\subsubsection{Possible inclusions}  It
will turn out that whenever $I, J$ are ideals, then
$\C_I\subseteq\C_J$ implies $I\subseteq J\subseteq \hat{I}$.
However, we will also see that this implication cannot be reversed
in general. Moreover, we will find ideal clones $\C_I, \C_J$ whose
unary fragments are comparable, i.e., $\C_I\cap\Oo\subseteq
\C_J\cap\Oo$, but which are incomparable, i.e., $\C_I\nsubseteq
\C_J$.

\subsubsection{Chains} We will prove that every ideal clone contains
a maximal subclone which is also an ideal clone. This implies that
there exist infinite descending chains of ideal clones. Using a
theorem from combinatorial set-theory, we will also show that there
exist infinite ascending chains of such clones.

\subsection{Overview}
Our paper is organized as follows: This introduction is followed by
a section on the mutual position of ideal clones, introducing the
important concept of $\hat{I}$ for an ideal~$I$
(Section~\ref{sec:mutualposition}). Precompleteness of ideal clones
and ways of testing it are  the topics of Sections~\ref{sec:maxreg}
and~\ref{sec:complexity}. In Section~\ref{sec:filters}, we move on
to the study of an alternative way of defining clones by means of an
ideal on~$X$, and how this new concept relates to the old one. All
this is done on countably infinite~$X$, but we then discuss the
possibility of a generalization to uncountable base sets in
Section~\ref{sec:uncountable}. Several problems that we had to leave
open are listed in Section~\ref{sec:openproblems}.

\subsection{Formal definition}
With the agenda at hand, we now officially introduce the main
objects of our interest.

\begin{defn}
If $I$ is an ideal on $X$ we denote
\[
\C_I:= \bigcup_ {n=1}^\infty\{f\in \O\un: f[A^n] \in I \mbox{ for all } A \in I\}
\]
and speak of the \emph{clone induced by} $I$. Clones of this form
will be called \emph{ideal clones}.
\end{defn}

\subsection{Support}\label{our.ideals}
 Let~$I$ be an ideal, and
 call the set of those elements of~$X$ which are contained in some
 set of~$I$ the \emph{support} of~$I$: $$\supp(I)= \bigcup I = \bigcup _{A\in I} A.$$
 For a subset
 $Y\subseteq X$, write $\pol(Y)$ for the set of all functions
 $f\in\O$ which satisfy $f(a_1,\ldots,a_n)\in Y$ for all
 $a_1,\ldots,a_n\in Y$. It is well-known and easy to see that
 $\pol(Y)$ is always a precomplete clone in $\Cl(X)$. Also, if we set
 $\proj(Y)$ to consist of all operations which behave like
 projections on~$Y$, then the interval
 $[\proj(Y),\pol(Y)]\subseteq\Cl(X)$ is isomorphic to $\Cl(Y)$ via
 the mapping $\sigma$ which sends every clone in the interval to the
 set of restrictions of its operations to~$Y$. Now setting
 $Y:=\supp(I)$, it is clear that $\C_I\in [\proj(Y),\pol(Y)]$, and
 that moreover $\sigma$ maps $\C_I$ to $\C_J$, where~$J$ is the
 restriction of~$I$ to~$Y$. Therefore, $\C_I$ can be imagined as an
 ideal clone on~$Y$, and consequently it is enough to understand
 clones induced by ideals which have \emph{full support}, i.e., which
 satisfy $\supp(I)=X$. This will be a permanent assumption from now
 on. Also, if~$I$ contains all subsets of~$X$, or only the finite
 subsets of~$X$, then $\C_I=\O$; therefore, we consider only ideals~$I$ satisfying $\P_{\rm fin}(X) \subsetneq I
\subsetneq \P(X)$, where $\P(X)$ denotes the power set of~$X$ and
$\P_{fin}(X)$ the set of all finite subsets of~$X$.

\subsection{Notation}

For a set of operations $\F\subseteq \O$, we denote the smallest
clone containing $\F$ by $\cl{\F}$. Given a clone $\C$ and a
function $f$ we let $\C(f)$ abbreviate $\cl{\C \cup\{f\}}$. If $n
\geq 1$, then $\F\un:=\F\cap\On$ is the set of all~$n$-ary
operations in $\F$. Given an operation $f\in\O$ and a subset
$A\subseteq X$, we will often write $f[A^n]$ for the image of the
appropriate power of~$A$ under~$f$, thereby implicitly assigning the
symbol~$n$ to the arity of~$f$. Rather than writing $\pi^1_1$, we
will use the symbol \emph{$id$} for the (unary) identity operation
on $X$.

\section{The mutual position of ideal clones}\label{sec:mutualposition}

     In this section we study order-theoretic properties of the mapping
     $I\mapsto \C_I$, and in particular examine the question of when
     $\C_I\subseteq \C_J$ holds for two distinct ideals $I,J$.

      We emphasize once again
     that throughout this paper, except for the special section on
     uncountable base sets, the base set $X$ is countably infinite.
     Moreover, all ideals are assumed to contain all finite sets, at
     least one infinite set, and do not contain all subsets of $X$.
     To stress these otherwise tacit assumptions we sometimes speak of
     `ideals in our sense'.

     Under these restrictions different ideals give rise to different
     clones. This follows from

     \begin{prop}\label{simpleconnections}
          For any ideals $I$ and $J$, $\C_I\uo\subseteq
         \C_J\uo$ implies $I\subseteq J$.
     \end{prop}

     \begin{proof}
         Assume $\C_I\uo\subset\C_J\uo$. In order to prove
         $I\subset J$, fix an infinite set
         $B\in J$ and let $A\in I$ be arbitrary.
         To show $A\in J$ choose a
         function $f:X\to A$ such that $f[B] = A$. As the range of $f$
         belongs to $I$, we have $f\in \C_I^{(1)}$. The assumption then
         yields $f \in \C_J\uo$, which in turn forces $A = f[B]$ into $J$.

 \end{proof}
     Next we show that $\C_I\uo \subseteq \C_J\uo$ is possible for
     distinct ideals, but only if $J\supseteq I $ is rather close to $I$.
     To make this precise we need the notion of a \emph{regular}
     ideal:

\begin{defn}
    Given an ideal~$I$, the \emph{regularization} of~$I$ is defined
    by $$\hat I := \{A\subset X: \mbox{ each infinite $B\subset A$
    contains an infinite set $C\in I$}\}.$$

   $I$ is called
   \begin{itemize}
   \item \emph{dense} (or \emph{tall})  iff $\hat I=\P(X)$, and
   \item \emph{regular}  (or \emph{nowhere tall}) iff $I=\hat I$.
   \end{itemize}
\end{defn}

\begin{rem} $\hat I$ can be written as $(I^\perp)^\perp$, where
$$ I^\perp = \{ A \subseteq X: \forall B\in I\; |A\cap B|<\aleph_0\}.$$

The initiated reader might notice that under \textsc{Stone} duality,
\begin{itemize}
\item ideals (in our sense) correspond to open sets in~$\beta\omega\setminus\omega$,
\item $\cdot^\perp$ gives the interior of the complement,
\item $I \mapsto \hat I$  corresponds to $U \mapsto \bar U^\circ$,
\item regular ideals just correspond to regular
open sets,
\item and similarly dense ideals to (topologically) dense open sets.
\end{itemize}

     Topological considerations will play no role in what follows but were
     of heuristic use in understanding the relationship between
     $I$ and $\hat I$.

\end{rem}

     The following facts are easy to verify and will be used without
     further reference.
     \begin{obss}\label{obs1}
      For all ideals $I$ we have  $I \subseteq \hat I$. If $I$ is not
     dense, then $\hat I$ is again an ideal (in our restricted sense)
     and turns out to be regular, i.e.\ $\hat{\hat I} = \hat I$. If $I
     \subseteq J \subseteq \hat I$, then $\hat J = \hat I$.
     Non-principal prime ideals are dense. The intersection of two
     dense ideals is dense.
      \end{obss}

     A little harder to prove is
     \begin{prop}
     All countably generated ideals are regular.
     \end{prop}

     \begin{proof}
      Let $I$ be a countably generated ideal. As $I \subseteq \hat I$
     is trivial, we have to prove $\hat I \subseteq I$.
     To do so, we consider an arbitrary $A\not\in I$ and find an
     infinite subset $B$ with no infinite $I$-subset.

     Let $(G_n)_{n=1}^\infty$ be an enumeration of some set of generators of $I$.
      $A\not\in I$ means that $A \setminus \bigcup_{k=1}^n G_k$
     is infinite for all $n$. So we can recursively build a sequence
     $(b_n)_{n=1}^\infty$ such that
     \[
     b_{n+1} \in A \setminus \left(\bigcup_{k=1}^n G_n \cup \{ b_1, \dots, b_{n}\}\right).
     \]
     Then $B:= \{b_n: n=1,2, \dots\}$ is the desired infinite subset
     of $A$, because by construction, no infinite subset of $B$ is
     covered by finitely many $G$s.
     \end{proof}

     \medskip

     \begin{prop}  \label{C_I_in_C_hatI} If $\C_I\uo \subseteq \C_J\uo$, then $J
     \subseteq \hat I$. If $I$ is not dense, then $\C_I \subseteq
     \C_{\hat I}$.
     \end{prop}

     If $I$ is dense, then $\hat I= \P(X)$ and
     the last implication is also true, except that
     we should not write $\C_{\hat I}$ in that case.

      \begin{proof}
     Assume $\C_I \uo \subseteq \C_J\uo$.         To see that $J\subset \hat I$, let $B$ in $J$ be arbitrary.
     We must show that below each infinite $C\subseteq B$ there is an
     infinite $I$-set. Assume not and let $C$ witness to this.
     Then any function $f$ mapping $C$ onto $X$ and being constant on $X
         \setminus C$ would belong to $\C_I$, but not to $\C_J$ (as
     $f[B] = X\not\in J$), a contradiction. Thus, there is no bad $C$
     and $B\in \hat I$, as demanded.

      To prove that $\C_I \subseteq \C_{\hat I}$, let $f\in \C_I$ and
     $A\in \hat I$. Our aim is to show that
         $f[A^n]\in \hat I$. Given an infinite subset $B$ of $f[A^n]$,
         pick $D\subset A^n$ such that $f[D]=B$ and such that $f$ is injective
         on $D$. If $\pi_1^n[D]\in I$, then set $A_1:=\pi_1^n[D]$ and $D_1:=D$.
         If $\pi_1^n[D] \in \hat I \setminus I$, then pick an infinite $I$-set $A_1\subset \pi_1^n[D]$ and let $D_1:= D \cap (A_1\times
         X^{n-1})$. In the second step thin out $D_1$ to an infinite set
         $D_2$ whose projection on the second coordinate lies in $I$.
         After $n$ such steps we arrive at an infinite set $D_n$ which
         is contained in $C^n$ for some $C\in I$. Thus $f[D_n]$ is an
         infinite $I$-set contained in $B$. Since $B$ was arbitrary, we
         conclude
         $f[A^n]\in \hat I$.
     \end{proof}

     \begin{ex}  Ideals $I,J$ showing that the implication $\C_I\uo \subseteq
     \C_J\uo \Rightarrow I \subseteq J \subseteq \hat I$ cannot be reversed.
     \end{ex}

     Let $X=\omega \times \omega$,
     \[
     I := \{ A \subseteq X:  \exists\ c\in\omega\ \forall\ n\in\omega\
      (\abs{A \cap (\omega\times\{n\})} \leq c)\}
     \]
     and
     \[
     J:=  \{ A \subseteq X:  \exists\ c,d\in \omega\ \forall\ n
     \in\omega\ \; (\abs{A \cap (\omega\times\{n\})} \leq c + d\cdot n)\}.
     \] Then
     $$\hat I=\hat J=\{A\subset X: \forall\ n\in\omega\ \;(A\cap (\omega\times\{n\}) \mbox{ is
     finite})\}.$$ Clearly $I\subset J\subset \hat I$. For every $1$-$1$ map
     $f:\omega\rightarrow \omega$ the map $F:X\rightarrow X$,
     $(k,n)\mapsto (k,f(n))$ preserves $I$, but it might not preserve $J$.
     Take for example
     $$ f(n)=\left\{\begin{array}{cl}
     2 m & \mbox{if $n=2^m$}\\
     2n+1 & \mbox{else.}
     \end{array}  \right.    $$

     The ideal $I$ of this  example
     was introduced in \cite{CH01} in order to show that there are
     non-precomplete ideal clones. As we shall see later in
     Proposition~\ref{ihatmax_no_choice},
     $\C_{\hat I}$ is the unique precomplete
     clone above $\C_I$. It is fairly easy to see that there are no
     other ideal clones between $\C_I$ and $\C_{\hat I}$.
     In
      fact, as was shown in \cite{GP08},  there is no clone at all
     between $\C_I$ and $\C_{\hat I}$, i.e.\ $(\C_I,\C_{\hat I})$ is a
     covering in the clone lattice.

     \begin{prop}\label{p8}
     Every ideal clone contains a maximal proper subclone which
     is itself an ideal clone.
     \end{prop}

     It follows that there are dense ideals whose clones are not precomplete.
     For, if $I$ is dense
     and  $\C_J$  a proper subclone of $\C_I$, then $X\in \hat I
\subseteq  \hat J$, by Proposition \ref{C_I_in_C_hatI}. So $J$
     is dense, but $\C_J$ cannot be precomplete.

     \medskip

This proposition will follow from  Lemmas~\ref{l7} and~\ref{l8}.  We
will need
the following notation:  For any ideal~$P$ on~$X$ and any function $f:X\to Y$
we define an ideal $\tilde f(P)$ on~$Y$ via $\tilde f(P)=\{B\subseteq Y:
f^{-1}[B]\in P\}$.   If $P$ is prime, then $\tilde f(P) $ is easily
     seen to be  a prime ideal again.
The quasiorder   $Q \le P \Leftrightarrow \exists f:   Q=\tilde
f(P)$ is called the {\sc Rudin-Keisler-}ordering on the prime ideals
on~$X$. It is, usually in the language of  ultrafilters, extensively
studied in the literature (see e.g.~\cite{CN74}).

\begin{lem}\label{l7}
    Let $\mathcal P$ be a set of non-principal prime ideals, and~$Q$ a
    non-principal prime ideal such that $\tilde f(P)\not=Q$ for all
    $P\in \mathcal P$ and all $f\in \Oo  $.  Denoting $I:=\bigcap
    \mathcal P$, we have $\C_{I\cap Q} \subseteq \C_I$.
\end{lem}

\begin{proof} Let $f\in \C_{I\cap Q}$ be~$n$-ary, $A\in I$.  We derive a contradiction from the assumption $B:=f[A^n]\notin I$.
So assume that $B\notin P$, for some $P\in \mathcal P$. Find a tuple
$\bar g = (g_1,\ldots, g_n)$ of unary functions such that $b = f( \bar g(b))$
for all $b\in B$, and $g_i[X]\subseteq A$ for $i=1,\ldots, n$.   From
$\tilde g_i(P)\not= Q$ we get $C_i\in Q$ such that
$g_i^{-1}[C_i]\notin P$.   Let $D:= B\cap g_1^{-1}[C_1]\cap \cdots
\cap g_n^{-1}[C_n]$, then $D\notin P$ (as~$P$ is a prime ideal).

On the other hand we have
$$
    g_i[D] \subseteq g_i[g_i\inv[C_i]]\subseteq C_i \in Q
    \mbox { and } g_i[D]\subseteq g_i[X] \subseteq A \in I.
$$
So $g_i[D]\in I\cap Q$, hence
$$
    D= \{ f(\bar g(d)): d\in D\}
    \subseteq f[g_1[D]\times \cdots \times  g_n[D]] \in I \cap Q,
$$
because~$f$ preserves $I\cap Q$.  But then $D\in P$, contradiction.

\end{proof}

\begin{lem}\label{l8}
Let $I \nsubseteq Q $ be  ideals,~$Q$ prime. \\If $f\in \C_I
\setminus \C_{I\cap Q}$, then $\C_I \subseteq  \C_{I\cap Q}(f)$.
\end{lem}

\begin{proof}
Assume first that the given $f\in \C_I \setminus \C_{I\cap Q}$  is unary.
Fix $A\in I\cap Q$
such that $B:= f[A]\notin I\cap Q$.   As $f\in \C_I$, we must
have $B\in I$, so $B\notin Q$, hence $X\setminus B \in Q$.
Note that~$B$ and~$A$ must be infinite.  By shrinking~$A$ we may additionally
assume
\begin{enumerate}
\item $f\restr A$  is a bijection from~$A$ onto~$B$.
\item $A\cap B = \emptyset$
\item There is some infinite $C\in I$ disjoint from $A\cup B$. \\
As $C \subseteq X\setminus B$, we have $C\in Q$.
\item Summarizing:~$A$,~$B$,~$C$ are disjoint,  $A, C\in I\cap Q $, $B\in I\setminus Q$.
\end{enumerate}
Now let $g:X\to X$ be constant outside~$C$ and map $C $ bijectively onto
$A\cup C$. Clearly $g\in \C_{I\cap Q}$. The function
$$
s(x,y,z):= \left\{\begin{array}{rl} y & \mbox{ if $x\in A$ }\\
z & \mbox{ if $x\in C$ }\\
x & \mbox{ otherwise  }
\end{array}\right.
$$
belongs to every ideal clone. It follows that
$p(x):= s(x, f(x), g(x))$ belongs to $ \C_{I\cap Q}(f)$.
Notice that~$p$ maps $X\setminus B$ bijectively onto~$X$; let~$q$ denote
its inverse.~$q$ is the identity outside $A\cup B \cup C$ and maps
$A\cup B \cup C $ to $A\cup C\in I$.   Hence $q\in \C_I$.

We now show that for arbitrary $g\in \C_I$ we have
$q\circ g\in \C_{I\cap Q}$:  Let $D\in I\cap Q$. As both
$q$ and~$g$ are in~$\C_I$, we have $q\circ g[D^m] \in I$.    But we
also have $q\circ g[D^m] \subseteq q[X] \subseteq X\setminus B\in  Q$.

Hence we have for all $g\in \C_I$:
$g= p\circ(q\circ g) \in  \C_{I\cap Q}(f)$,
so $\C_I \subseteq   \C_{I\cap Q}(f)$,
as desired.

Now consider the case that
the given $f  \in \C_I \setminus \C_{I \cap Q}$ is $n$-ary for
some $n>1$.
We show that $f$ can be replaced by a unary function $f'$.
 For
some $A\in I \cap Q$ we must have $f[A^n] \not\in I \cap Q$. Choose
$g_1, \dots, g_n : X \to A$ such that $\{(g_1(a), \dots, g_n(a)):
a\in A\} = A^n$. Having range $A$, all $g_i$ belong to $\C_I$ and to
$\C_{I\cap Q}$. It follows that the unary function $f' = f(g_1,
\dots, g_n)$ belongs to $\C_I$ and to $\C_{I\cap Q}(f)$ but not to
$\C_{I \cap Q}$ because $f'[A] = f[A^n] \not\in I \cap Q$.  From the
first case it follows that $\C_I \subseteq \C_{I\cap Q}(f')
\subseteq \C_{I\cap Q}(f)$.
\end{proof}

  \begin{proof}[Proof  of Proposition~\ref{p8} from the lemmas]
     Let $I$ be the given ideal. It is well known that $I$ is an
     intersection of prime ideals. Therefore, we  can choose for each $A\not\in I$
     some prime ideal $P_A$ such that $I \subseteq P_A\not\ni A$.
     Then $I = \bigcap_{A\not\in I} P_A$ and the set ${\mathcal P}:= \{P_A: A\not\in
     I\}$ has  power at most $2^{\aleph_0}$. As there are only
     $2^{\aleph_0}$ functions $X\to X$, the set $\{\tilde g(P): P \in
     {\mathcal P} \mbox{ and } g \in \O^{(1)}\}$ has power $2^{\aleph_0}$,
     too.

     As $I$ is an ideal `in our sense', there exists some infinite
     $B\in I$ such that $C:= X \setminus B$ is also infinite. Again
     it is well-known that there are $2^{2^{\aleph_0}}$ non-principal
     prime ideals containing $C$. So letting $Q$ be one of them,
     which is not of the form $\tilde g(P)$ for any $P\in \mathcal P$, we
     have $I \not\subseteq Q$ (because of $B$) and the two lemmas now
     show that $J:= I \cap Q$ is as desired.
\end{proof}

     \begin{cor}
     Below each ideal clone there is an infinite descending chain of
     ideal clones. There are also infinite ascending chains of ideal clones.
     \end{cor}
     We do not know what other types of chains of ideal clones exist.
     \medskip

     To get the descending chain, start with an arbitrary ideal clone
     and successively apply  Proposition~\ref{p8}.

     To get the ascending chain is much harder because a deeper
     result of combinatorial set-theory is needed:
\begin{thm}[{\textsc{Kunen} \cite{Kun72},  see also \cite[Theorem 10.4]{CN74}}]
There are $2^{\aleph_0} $ many maximal ideals
on $\N$ which are pairwise incomparable in the {\sc Rudin-Keisler-}order.
\end{thm}

We only need countably many incomparable prime ideals.  So let
     $(P_k)_{k=1}^\infty$ be a
     sequence of prime ideals
     such that for any
     $g\in \O^{(1)}$ and any $m \not= n$ we have $\tilde g(P_m) \not=
     P_n$.

     Now we may simply put $I_n:=
     \bigcap_{k\geq n} P_k$. Then $I_n = I_{n+1}\cap P_n$ and Lemma~\ref{l7} says $\C_{I_n} \subsetneq \C_{I_{n+1}}$.

\leftskip=0cm

\section{Precompleteness and regularity}\label{sec:maxreg}

Problem~1 of \cite{CH01} asks whether every ideal clone lies below a
precomplete ideal clone. The following theorem gives a positive
answer to this question and shows that an even stronger statement
holds.

\begin{thm}\label{solution2question1}
     Every ideal clone is contained in a precomplete clone. Moreover,
    every precomplete clone containing an ideal clone is an ideal
    clone itself.
\end{thm}

To prove Theorem~\ref{solution2question1} we
 employ the following lemmas which will also be
useful later on.

\begin{lem}\label{completenesslemma}
Let $I$ be an ideal and $f\in \O$. If $f[B^n]=X$ for some $B\in I$,
then $\C_I(f)= \O$.
\end{lem}

\begin{proof}
    Choose $g_1,\ldots, g_n: X \to B$ such that $(g_1,\ldots, g_n)
    :X\to X^n $ maps $B$ onto $B^n$.
    Then the function $h:= f(g_1, \ldots, g_n)$ maps $B$ onto $X$.
    Pick  $C\subset B$ such that $h$ maps $C$ bijectively onto $X$
    and  fix  functions $h_k\in\O^{(k)}$ mapping $X^k$
    bijectively onto $C$.
    Let $g$ be any $k$-ary operation on $X$. Its action on
    $\bar x$ can be channeled through $C$ as follows: For every
    $\bar x$, there exists precisely one $d\in C$ such that
    $h(d)=g(\bar x)$. Moreover, there is precisely one $c\in C$ such
    that $h_k(\bar x)=c$. We define a unary operation $\tilde g: C \To
    C$ by setting $\tilde g(c)=d$ for all $c,d\in C$ obtained that
    way. Now extend $\tilde g$ anyhow to a unary operation on $X$, but in such a way that its range is still contained in $C$.
    We have: $\bar x \mapsto h_k(\bar x) =: c  \mapsto \tilde g(c):= d \mapsto
    h(d) = g(\bar x)$.
    From $g(\bar x) = h(\tilde g(h_k(\bar x)))$ we conclude that
    $g$ is a composition of $f$ and $g_1, \dots, g_n, h_k, \tilde g$. As
    the latter functions all have ranges in $I$, namely $B$ or $C$, they
    automatically belong to $\C_I$. So $g\in \C_I(f)$. As $g$ was
    arbitrary, $\C_I(f) = \O$, as claimed.
\end{proof}

\begin{lem}\label{constructionlemma}
 Let $A\subseteq X$ be  an infinite set and $\A$ a clone containing an ideal clone.
 Set $J:=\{B\subseteq f[A^n]:f\in \A\}$. Then either $X\in J$,  or:
  $J$ is
 an ideal (in our sense), $A\in J$,   and $\A\subseteq \C_J.$
\end{lem}

\begin{proof}
We start by showing that $J$ is an ideal: By definition $J$ is
closed under the formation of subsets.
 To see that $J$ is closed under finite
unions, consider the switching function
$$
    s(x,y,u,v) =
    \left\{\begin{array}{rl} x,& u=v\\y,& u\not=v.\end{array} \right.
$$
One easily checks that every ideal clone contains $s$ and thus $s\in
\A$. Assume that $B,C\in J$, say $B\subset f[A^m], C\subset g[A^n]$
for some $f,g\in \A$. Then $B\cup C \subset f[A^m]\cup
g[A^n]=s[f[A^m]\times g[A^n]\times A\times A]$. Therefore, $B \cup C
\in J$ as the rightmost set is the image of $A$ under the operation
$s(f(x_1,\ldots,x_m),g(y_1,\ldots,y_n),u,v)\in\A$. Since every ideal
clone contains all functions with finite image, so does $\A$.
Therefore $J$ contains all finite sets. As $A = \id[A] \in J$ there
is an infinite set in $J$. So $J$ is an ideal unless $X\in J$ (in
which case $J = \P(X)$).

Finally we have to prove that $\A\subset \C_J$. Consider any  $f\in \A^{(n)}$
and $B\in J$. Choose $g\in \A$ such that $B\subset g[A^m].$ Then
$f[B^n]\subset f[g[A^m]\times\ldots \times g[A^m]]\in J$ since the
latter set is the image of $A$ under the operation
$f(g(x_1^1,\ldots,x_m^1),\ldots,g(x_1^n,\ldots,x_m^n))\in\A$.
\end{proof}

\begin{proof}[Proof of Theorem~\ref{solution2question1}]
    Pick an infinite $A\in I$ and choose a function $g\in\O\uo$ which maps
    $A$ onto $X$ and is constant on $X\setminus A$.
    By Lemma~\ref{completenesslemma}, $\O= \C_I(g)$.
    Consider the set $S$ of all clones above $\C_I$ which do not contain
    $g$.
    {\sc Zorn}'s Lemma yields the existence of a maximal element
    $\A$ in~$S$. Every clone which is strictly larger than~$\A$ contains~$g$
    and is thus equal to~$\O$. Hence $\A$ is precomplete.

    It remains to show that each precomplete $\A$
    above $\C_I$ is an ideal clone.
    Set $J=\{B\subset f[A^{n}]: f\in \A\}$.

    $X$ cannot be in $J$, for, otherwise  $X=f[A^{n}]$ for some $f\in \A$. By Lemma~\ref{completenesslemma} this implies
    $\O=\C_I(f) \subset \A $ which is impossible
    since $\A$ is a proper subclone of $\O$.
    So, by Lemma~\ref{constructionlemma}, $J$
    is an ideal and $\A \subseteq \C_J \not=\O$
     Since $\A$ is precomplete, we must have
    $\A=\C_J$.
\end{proof}

The following proposition reveals the significance of regular ideals
for our purposes.

\begin{prop}\label{ihatmax_no_choice}
    Assume that $I$ is not dense. Then $\C_{\hat I}$ is
precomplete and in fact the only precomplete clone above $\C_I$.
So $\C_I$ is precomplete iff $I$ is regular.
\end{prop}

 In the case of non-dense ideals we therefore have a
 very satisfactory criterion:  the
precompleteness of $\C_I$ can be read off the ideal $I$ directly
without even looking at the operations in $\C_I$. If $I$ is dense,
then no general statement can be made:  Prime ideals are
obviously dense
and give rise to
 precomplete clones (Example~\ref{primecomplete}).  Proposition~\ref{p8}
 and the remark following it show that there are many dense ideals
which lead to a
non-precomplete clone.

Notice that the precompleteness of $\C_{\hat I}$ is proved in ZF, i.e.\ without the use of the
Axiom of Choice.

\begin{proof}
    From $X\not\in \hat I$ we infer that $\hat I$ is an ideal (in
    our sense) to which the former machinery applies.
    We start by proving the precompleteness of $\C_{\hat I}$.
    Consider any $f\in \O^{(n)} \setminus
    \C_{\hat I}$. Pick $B\in \hat I$ such that $f[B^n]$ contains an
    infinite set $C$ which contains no infinite $I$-set. Choose any
    $g\in\O\uo$ such that $g[C]=X$ and such that $g$ is constant on $X\setminus
    C$. Then  $g\in \C_I\subseteq \C_{\hat I}$ (as images of
    $I$-sets are finite) and $g\circ f$ maps $B^n$ onto $X$.
    Thus $\C_{\hat I}(f) \supseteq  \C_{\hat I}(g\circ f)=\O$, by
    Lemma~\ref{completenesslemma}.

    Next we consider any precomplete
    $\A\supseteq \C_I$
    and show that $\A = \C_{\hat I}$. From Theorem~\ref{solution2question1}
    we know that $\A = \C_J$ for
    some ideal $J$. But then $I \subseteq J \subseteq \hat I$, by
    Propositions~\ref{simpleconnections} and~\ref{C_I_in_C_hatI}.
    Thus, $\hat J = \hat I$ (confer Observations~\ref{obs1}). So, $\A = \C_J
    \subseteq \C_{\hat J} = \C_{\hat I}$. By maximality, $\A =
    \C_{\hat I}$, as desired.
\end{proof}

It was already  mentioned that  there are $2^{2^{\aleph_0}}$
precomplete clones on a countable base set. The hitherto known
constructions of that many clones all used the Axiom of Choice in
one way or the other. Next we set out to produce the maximal number
of precomplete clones without using AC.

\begin{thm}[ZF]\label{many_clones_no_choice}
    There exists an injective mapping from $2^{2^\omega}$ into the set of precomplete clones over a countable set.
\end{thm}

 From the above, it is clear that we only have to construct many
regular ideals.

We begin with the following \begin{lem} Let $\cal B$ be a
collection of subsets of $X$ and denote by $I_{\cal B}$
the ideal generated by $\cal B$ and all finite sets. If $A$ is
infinite and $A\in \hat I_{\cal B}$, then $A$ has infinite
intersection with some $B\in\cal B$.\end{lem}

\begin{proof} By definition, there is some infinite $C\subseteq
A$, such that $C\in I_{\cal B}$, i.e.\ $C \subseteq F \cup B_1 \cup
\dots \cup B_k$ for some finite set $F$ and finitely many $B_1,
\dots, B_k\in {\cal B}$. As $C$ is infinite, $C\cap B_i$ must be
infinite for one $B_i$. Then $A\cap B_i$ is all the more infinite.
\end{proof}

\medskip

For the actual construction we consider an infinite  family $\A$ of what are
called almost disjoint subsets of $X$, i.e., the sets in $\A$ are
infinite but have finite pairwise intersections.
If $\cal B$ is a  subcollection of $\cal A$,
the above lemma shows $\hat I_{\cal B} \cap \A = {\cal B}$. It
follows that different ${\cal B}$s give rise to different regular
ideals, hence to different precomplete clones (in fact for
${\cal B} = \emptyset$ or ${\cal B} = \A$ we may not get ideals
in our sense).
In other words, from $\A$ we have constructed $2^{|\A|}$
distinct precomplete clones.

We are left with constructing an almost disjoint family of power
$2^{\aleph_0}$ without using the Axiom of Choice. There are a number
of such constructions in the textbooks of set-theory. The most
popular one (due to {\sc Sierpi\'nski}) is based on  $X:= 2^{<
\omega}$, the countable set  of all finite 0-1-sequences. If $\xi$
is an infinite 0-1-sequence, let $A_\xi$ be the set of its initial
segments. Then
 $\{A_\xi : \xi \in 2^{\omega}\}$ is an almost
disjoint family. Clearly, each infinite sequence has infinitely many
initial segments. So each $A_\xi$ is infinite. Moreover, if $\xi
\not=\eta$ then $\xi(n) \not= \eta(n)$  for some (first) $n$. Then
$A_\xi$ and $ A_\eta $ have no initial segments of length $\geq n$
in common. Hence $A_\xi \cap A_\eta$ is finite.

\section{The complexity of testing precompleteness}\label{sec:complexity}

The following precompleteness test was already proved in
\cite{CH01}:
\begin{quote}{\it
    $\C_I$ is precomplete iff for all $A\not\in I$ there exists $f\in
    \C_I$ such that $f[A^n]=X$.}
\end{quote}
With the tools from the previous section the proof becomes  very
short. If $\C_I$ is not precomplete, then $\C_I \subseteq \C_J$ for
some ideal $J \supsetneq I$. Taking $A\in J\setminus I$ and any
$f\in \C_I$ we cannot have $X= f[A^n]$, for,   $f$ would belong to
$\C_J$, too, and $A\in J$ implies $f[A^n] \in J\not\ni X$.

In the other direction we let $A\not\in I$ be as stated.
Applying Lemma~\ref{constructionlemma} to $\C_I$ and $A$ we get
$J \supset I$. The assumption on $A$ just says $X\not\in J$. So
$J$ is an ideal and
 $\C_J \supsetneq \C_I$, disproving maximality of  $\C_I$.

\medskip

Applications of the test can be simplified by the observation that
instead of $f[A^n] = X$ it is sufficient to demand that $f[A^n]$ be
big, i.e.\ has a complement in $I$, or, more formally, belong to
$F:= \{X \setminus B: B\in I\}$, the dual filter of $I$. This is
justified by the following

\begin{lem}\label{lem:mapOnesetToX}
 Each set in $F$ can be mapped onto the whole of $X$  by a
unary $\C_I$-function.
\end{lem}

\begin{proof} Denote the set in question by $B$. Belonging to  $
F$, $B$ must be infinite.
If $B$ does not contain
any infinite subset belonging to $I$, then any function mapping
$B$ onto $X$ and constant on $X\setminus B$ will do, because it
is automatically in $\C_I$. Otherwise, there is an infinite
$C\subseteq B$ such that $C\in I$. Let us split it: $C= C_1\cup
C_2$ with both parts infinite. Then any function $f$ such that
$f[C_1] = C$, $f[C_2] = X \setminus B$, and $f(x) = x$ on $X
\setminus C$  maps $B= (B \setminus C) \cup
C_1\cup C_2$ to $(B\setminus C) \cup C \cup (X\setminus B) = X$.
Any such $f$ also belongs to $\C_I$, because from $D\in I$ we have
\[
f[D] = f[C_1 \cap D ] \cup f[C_2 \cap D] \cup f[D\setminus C]
\subseteq C \cup (X\setminus B) \cup (D\setminus C) \in I.
\]
\end{proof}

\begin{ex}\label{primecomplete}
Taking this simplification into account, precompleteness for $\C_I$
with $I$ prime becomes a triviality: If $A\not\in I$, then
$\id[A]$ is in the dual (ultra)filter.
\end{ex}

There is another way to liberalize the condition in the test.
Instead of one $\C_I$-function one can allow a finite (and with some
care even an infinite) number of functions. In the rest of this
section we shall be mainly concerned with the question of how many
functions of which arities are needed to test the precompleteness of
a given $\C_I$. It will turn out that binary functions do not
suffice, as some of us hoped to prove before we knew the
counterexample. To make things more precise, we introduce some
notation. For $A\subseteq X, n\geq 1$ and $1 \leq p \leq \infty$ we
write $T(A,n,p)$ iff there is a sequence $(f_k)_{k=1}^p$ of $n$-ary
functions such that

\[
(i)\; \bigcup_{k=1}^p f_k[A^n] = X \qquad \mbox{ and } \qquad
(ii)\;  \bigcup_{k=1}^p f_k[B^n] \in I \mbox{ for all } B\in I.
\]

Obviously, $T(A,n,p)$ also depends on $I$. But the ideal will
usually be fixed in the context. Therefore, we can safely
suppress it in the notation.

If $p$ is finite, then $(ii)$ just says  that all $f_k$
belong to $\C_I$.
By the above lemma, condition $(i)$ can be replaced by
$\bigcup_{k=1}^p f_k[A^n] \in F$.

We now let $p$ enter the game and call the result

\medskip

\begin{thm} \label{MaxCri}
$\C_I$ is precomplete iff for each $A\not\in I$
there are $n \geq 1$ and $1\leq p \leq \infty$ such that
$T(A,n,p)$ holds.
\end{thm}

\begin{rem}
This criterion strongly depends on our assumption that the base set
is countable.
Example~\ref{prop:uncountable:boundedIdealPrecompleteClone} in
Section~\ref{sec:uncountable} shows that this criterion fails on all
uncountable sets.

\end{rem}

In all examples we know,  a stronger form holds: there are $n$ and
$p$ such that $T(A,n,p)$  for all $A\not \in I$.
It is an open problem if that is always the case. If it is, we
say in a somewhat sloppy way,  that  $\C_I$ is precomplete via
$p$ $n$-ary functions.

\medskip
\begin{proof}[Proof of Theorem~\ref{MaxCri}]
Using the new notation the precompleteness test mentioned at the
beginning of this section
 reads: $\C_I$ is precomplete iff for all   $A \not\in I$ there is some $n$ such that
 $T(A,n,1)$ holds.
The rest follows from the  observation that for all infinite $A$,
$n\geq 1$ and $1< p < q <\infty$
\[
T(A,n,1) \Rightarrow T(A,n,p) \Rightarrow T(A,n,q) \Rightarrow T(A,n,\infty)
\Rightarrow T(A,n+1,1).
\]
Below we give examples showing that none of these implication can be reversed.

The others being obvious, only the last implication needs proof.
Let $(f_k)_{k=1}^\infty$ witness $T(A,n,\infty)$. As $A$ is
infinite, we can choose a sequence $(a_k)_{k=1}^\infty$ of
pairwise distinct elements of $A$. It allows us to define
\[
f(x_1, \dots, x_n,y):= \left\{\begin{array}{rl} f_k(x_1, \dots,
x_n),& y= a_k\\
y, & \mbox{otherwise.}
\end{array}\right.
\]
An easy verification then shows that this function witnesses $T(A,n+1,1)$.
\end{proof}

\medskip

 \cite[Problem D]{GP06survey} asks for a precompleteness test using unary
functions only. The following proposition gives such a test,
which seems however of little practical use. The case $n=1$ will
turn out to be important in the next section, though.

\begin{prop}\label{unarytest}
 $\C_I$ is precomplete iff for each $A\not\in I$ there exist
$n\geq 1$ and unary functions $g_1, \dots, g_n : X \to A$ such
that $g_1^{-1}[B] \cap \dots \cap g_n^{-1}[B] \in I$ for all
$B\in I$.
\end{prop}

The proposition will follow from  Theorem~\ref{MaxCri} and the following
lemma,
by  assembling  the $g_i$ to a vector function $\bar g= ( g_1,
\dots, g_n)$.

\begin{lem}\label{preim}
 If $A\subset X$ is infinite, then condition $T(A,n,\infty)$
is equivalent to the existence of a
 vector
function $\bar g: X \to A^n$ such that $\bar g^{-1}[B^n]\in I$
for all $B\in I$.
\end{lem}

\medskip

\begin{proof} Let $\bar g: X \to A^n$ be as stated. Choose $c\in X$ arbitrarily.
Define a sequence $(f_k)_{k=1}^\infty$ of $n$-ary functions in
such a way that
\[
\{ f_k(\bar y): k=1,2,\dots \} = \{c\} \cup \{x\in X:
\bar g (x) = \bar y\}
\]
for each tuple $\bar y\in X^n$. This is possible because the set of
the right-hand side is  non-empty and countable.

These functions are as demanded by $T(A,n,\infty)$. Indeed, for each
$x\in X$, the tuple $\bar g(x)$ belongs to $A^n$ and has $x$ in its
preimage. So $x = f_k(\bar g(x)) \in f_k[A^n]$ for some $k$. Hence
$X\subseteq \bigcup_{k=1}^\infty f_k[A^n]$. Moreover, for any $B\in
I$ we have $\bigcup_{k=1}^\infty f_k[B^n] \subseteq \{c\} \cup \bar
g^{-1}[B^n]\in I$, as demanded.

In the other direction we start from a sequence $(f_k)_{k=1}^\infty$
as in $T(A,n,\infty)$. As $\bigcup_{k=1}^\infty f_k[A^n]= X$, for
each $x\in X$ there exists some (smallest) $k$ and some $\bar a\in
A^n$ such that $f_k(\bar a) = x$. Let $\bar g(x)$ pick such tuple
$\bar a$. Then $\bar g$ is a vector function $X \to A^n$ and for
each $B\in I$ we have  $\bar g^{-1}[B^n] \subseteq
\bigcup_{k=1}^\infty f_k[B^n] \in I$.
\end{proof}

In the rest of this section we give examples showing that there
is no general bound on the number and arity of the functions
that are needed to establish precompleteness according to
Theorem~\ref{MaxCri}.

\begin{ex}\label{prop:criterion:manyUnaryOperations}
    For every $1<p\leq \infty$, there exists an ideal $I$ such that
$\C_I$ is precomplete via $p$ unary functions but
fewer unary functions do not suffice.
\end{ex}

Let $Y$ be any countable set and put $X:= \bigcup_{k=1}^p Y
\times \{k\}$. Let $P$ be any prime ideal on $Y$. Then the ideal
$I$ we aim at, consists of all sets $A= \bigcup_{k=1}^p A_k \times
\{k\}\subseteq X$ for which all $A_k$ belong to $P$.

First we show that $T(A,1,p)$ holds for all $A\not\in I$. Write
$A = \bigcup_{k=1}^p A_k\times \{k\}$ and notice that for some
$k$ we must have $A_k \not\in P$. For notational convenience we
assume that $A_1\not\in P$. By Lemma~\ref{lem:mapOnesetToX}
there is a $\C_P$-function $f:Y\to Y$ that maps $A_1$ to $Y$
(not being in the prime ideal $P$ means being in its dual (ultra)filter).

Using this $f$ we define $f_k: X \to X$ by setting
\[
f_k(y,m) = \left\{\begin{array}{rl} \left( f(y),k\right),& m=1\\
(y,m),& m\not=1\end{array}\right..
\]
Then
\[
\bigcup_{k=1}^p f_k[A] \supseteq  \bigcup_{k=1}^p f_k[A_1\times \{1\}]
= \bigcup_{k=1}^p f[A_1] \times \{k\}  = \bigcup_{k=1}^p Y \times \{k\} = X
\]
and  for any $B = \bigcup_{m=1}^p B_m\times \{m\}\in I$
\[
\bigcup_{k=1}^p f_k[B]
\subseteq B \cup \bigcup_{k=1}^p f_k[B_1\times \{1\}] =
B \cup \bigcup_{k=1}^p f[B_1]\times \{k\}
\in I.
\]

Next we show that $T(Y \times \{1\}, 1, q)$ does not hold for
any $q<p$ (notice that $q$ is finite, even if $p = \infty$).
Consider functions $g_1, \dots, g_q: X \to X$ such that
$\bigcup_{m=1}^q g_m[Y \times \{1\}] = X$. We show that
$g_m\not\in \C_I$ for  at least one of the functions.
\medskip

For each $m\leq q$ we write $g_m[Y \times \{1\}]$ as
$\bigcup_{k=1}^p D_{mk}\times \{k\}$. Then for each $k$ we have
$\bigcup_{m=1}^q D_{mk} = Y$.  It follows that for each $k$ there must be some
$m$ such that $D_{mk} \not\in P$. As there are more $k$s than
$m$s, one $m$ must serve two $k$s. In other words: there is some
$m\leq q$ such that $D_{mi}\not\in P$ and $D_{mj}\not \in P$ for some $i\not=j$.

The sets $D_{mi}\times \{i\}$ and $D_{mj}\times \{j\}$ do not belong
to $I$ then. As they are disjoint, so are their preimages under
$g_m$ and even more so
\[
g_m^{-1}[ D_{mi}\times \{i\}] \cap \left( Y \times \{1\}\right)
\quad \mbox{ and }\quad
g_m^{-1}[ D_{mj}\times \{j\}] \cap \left( Y \times \{1\}\right)
\]
As (the trace of) $I$ is prime on $Y\times \{1\}$ one of these
disjoint sets is in $I$, but none of their $g_m$ images is. So $g_m
\not\in \C_I$, as claimed.

As a byproduct of the above, we get
\begin{ex}\label{unary_part_smaller_binary_not}
Ideals $I,J$ such that $\C_I\uo\subseteq
    \C_J\uo$ but      $\C_I\not \subseteq \C_J$.
    \end{ex}

    Let  $I$ be the ideal constructed in the previous example
for $p= \infty$ and let  $B\not\in I$ be such that $T(B,1,q)$
fails for all finite $q$.
Then $\C_I$ is precomplete and
    $J=\{A\subset \bigcup_{i=1}^n f_i[B]:f_1,\ldots,f_n\in \C_I\uo\}$
    is an  ideal in our sense.
The definition of $J$ immediately yields   $\C_I\uo\subseteq
    \C_J\uo$ and, since $B\in J\setminus I$, this inclusion is proper, by Proposition~\ref{simpleconnections}.
    The precompleteness  of $\C_I$ shows that we cannot have
    $\C_I\subseteq \C_J$.

\begin{ex}\label{maximal_but_not_via_unary}
For each $n\geq 1$ there is an  ideal $I$ such
that  $\C_I$ is precomplete via one $n+1$-ary function but not
via (even infinitely many) $n$-ary functions.
\end{ex}

We let $P$ be a non-principal prime ideal on the countably
infinite set $Y$ and put
\[
X := Y^{n+1} \quad \mbox{ and }\quad I = P^{n+1} = \{ B
\subseteq X: B \subseteq C^{n+1} \mbox{ for some } C \in P\}.
\]
Regardless of the choice of $P$ we then have $T(A,n+1,1)$ for all
$A\not\in I$.
To see this, let $A\not\in I$ be given. By definition,  there is a
projection $p: X= Y^{n+1} \to Y$ on one of the coordinate axes such
that $B:= p[A]$ does not belong to $P$.
By Lemma~\ref{primecomplete}, we may choose a function $g: Y \to Y$
that maps $P$-sets to~$P$-sets and $B$ onto the whole of~$Y$.
Then
\[
f(x_1, \dots, x_{n+1}) = \Big(\, g(p(x_1)), \dots, g(p(x_{n+1}))\,\Big)
\]
belongs to $\C_I$ and maps $A^{n+1}$ to $X= Y^{n+1}$.

We do not know if there are prime ideals $P$ such that the above
construction would yield an ideal $I$ such that $\C_I$ was
precomplete via infinitely many $n$-ary functions. By properly
choosing   $P$, we show that, at least, this is not always the case.

We first explain what property of $P$ yields the desired result. We
let $0$ denote any fixed element of $Y$ and consider  $ A:= \{(0,0,
\dots, 0, y): y\in Y\}.$ Then $A$ has full projection onto the last
coordinate, hence $A\not\in I$.

If $T(A,n,\infty)$ were true, then Lemma~\ref{preim} would yield a
vector function $\bar g = (g_1, \dots, g_n): X \to A^n$ such that $\bar g^{-1}[B^n]
\in I$ for all $B\in I$.
Let $p$ denote the projection $Y^{n+1} \to Y$ onto
the last coordinate (so $p[A] = Y$) and put $h_i = p(g_i)$. Then
$\bar h = (h_1, \dots, h_n)$ maps $X$ to $Y^n$. Moreover,
for each
$B\in P$ we have $\tilde B:= \{(0,0,\dots, 0,b): b\in B\} \in I$
and, therefore,
\[
\bar h^{-1}[B^n] = \bigcap_{i=1}^n h_i^{-1}[B] =
\bigcap_{i=1}^n g_i^{-1}[\tilde B] = \bar g^{-1}[\tilde B^n] \in I.
\]

So, in order to guarantee that $T(A,n,\infty)$ fails, we
 can choose $P$ in such a way that
\BE
\xitem{*} for each vector function $\bar h: X = Y^{n+1}
\to Y^n$ there exists some $C\in P$ such that
$\bar h^{-1}[C^n] = h_1^{-1}[C] \cap \dots \cap h^{-1}_n[C]$
does not belong to $I$ (e.g. because this set has full projection
 onto some coordinate).
\EE

Let $H$ denote the set of all $\bar h: Y^{n+1} \to Y^n$ to be
dealt with.  For each $\bar h\in H$
we are going to choose
a set $C_{\bar h}$ such that $\bar h^{-1}[{C_{\bar h}}^n]$ has full
projection onto one coordinate. Some extra care is needed
to ensure that all the chosen sets peacefully live together in
one prime ideal.

To achieve this, we need large independent families: A family
$\Bbb{F}$ of subsets of $Y$ is called \emph{independent} iff for all
finite disjoint $\Bbb{F}_1,\Bbb{F}_2\subseteq \Bbb{F}$
$$
    \bigcap\{F: F\in\Bbb{F}_1\}\, \cap\, \bigcap \{Y\setminus F:F\in \Bbb{F}_2\}\neq\emptyset.
$$
The following classical result shows that large independent families
exist:
\begin{lem}[{{\sc Fichtenholz, Kantorovich and Hausdorff}, see \cite[Lemma 7.7]{Jec02}}]
There is an independent
family $\Bbb{F}$ of subsets of $Y$ which has power $2^{\aleph_0}$.
\end{lem}
A particularly concrete example of such a
family was given by {\sc Menachem Kojman} in \cite{GK}:
on the countable set $\Bbb Z[x]$
of polynomials with integer coefficients  define, for
any real $r$, the set $Z_r$ as
$$Z_r:= \{ p(x)\in \Bbb Z[x]: p(r)>0 \}. $$
Then  the collection of all $Z_r$  with $r \in \Bbb R$ is independent.

As $|\Bbb{F}| = |H|= 2^{\aleph_0}$ we can split   $\Bbb{F}$ as
$\bigcup_{\bar h \in H} \Bbb{F}_{\bar h}$, where the
$\Bbb{F}_{\bar h}$ are pairwise disjoint and  infinite.
If we choose   the $C_{\bar h}$ as Boolean combinations
of sets in $\Bbb{F}_{\bar h}$ (and different from $X$), then  no finite union of them covers $Y$. In
other words, $\{C_{\bar h} : \bar h \in H\}$ extends to a
prime ideal.

It remains to explain, how an individual $C_{\bar h}$ can be found.
So let
 $\bar h = (h_1, \dots, h_n) : Y^{n+1} \to Y^n$ be given   and take a
decomposition  $Y = D_1 \cup  D_2 \dots
\cup D_{n+1}$ into mutually disjoint non-empty sets (Boolean
combinations from sets in $\Bbb{F}_{\bar h}$).
We claim that
for one index $m \leq n+1$ the preimage  $\bar h^{-1}[(Y\setminus D_m)^n]\subseteq Y^{n+1}$ has full
projection onto one coordinate.
Then the corresponding set $C_m:= Y \setminus D_m$ does the required job.

Consider
any $x\in X$. By the Pigeonhole Principle,
there must be some $i\leq n+1$ such that
$D_i$ contains none of $h_1(x), \dots h_n(x)$, so $C_i$ contains
all of them, or,  in other words,
$\bar h(x) \in C_i^n$.
This shows
\[
Y^{n+1} = X = \bigcup_{i=1}^{n+1} \bar h^{-1}[ C_i^n].
\]
But then one set in the union has full projection. For, if for
each $i$ there were some $y_i$ that does not occur as $i$-th
coordinate of an element of $\bar h^{-1}[C_i^n]$, then $(y_1, y_2,
\dots, y_{n+1})$ would not occur in the whole union.

\section{The dual construction: clones from filters}\label{sec:filters}
There are other ways to construct clones from ideals. In this
section we discuss one approach  which is, in a sense, dual to the
former one.
Functions have been put in $\C_I$ if images of small sets are small.
Now we consider those functions for which preimages of small
sets are small.
Taken literally, this idea would lead to the set
\[
\bigcup _{n\geq 1} \{ f\in \O^{(n)}: f^{-1}[A] \in I^n \mbox{ for all }
A\in I\},
\]
which is not a clone, however (e.g. because it does not include the
projections). The generated clone, denote it by $\D_I$, consists of
all functions that are essentially in the set, i.e.\ up to
fictitious variables. We did not study this construction in detail
because there is another more promising way to make the above idea
precise.
 The first step is to pass to complements:
preimages of big sets have to be big.
In the context of a given ideal, the small  subsets of $X$
are those in $I$ and the big ones those in the dual filter
$F:=\{X\setminus A: A\in I\}$.
If we take subsets of $X^n$ to be big if they belong to
$F^n$ (the filter generated by all $B^n$ for $B\in F$), then we
arrive at the following set of functions
\[
\S_F := \bigcup_{n=1}^\infty \{ f\in \O^{(n)}: f^{-1}[B]\in F^n  \mbox{ for all } B
\in F\}.
\]
A straight-forward verification shows $\S_F$ to be a clone.
The functions in $\S_F$ will be called $F$-continuous. An
 alternative description of $F$-continuity demands that  for all
$B\in F$ some $C\in F$ can be found such that $f[C^n] \subseteq B$.

An easy verification shows that $\D_I \subseteq \S_F$ while
$\D_I^{(1)} = \S_F^{(1)}$. Later on we have several times occasion
to test wether a given unary function belongs to $\S\uo_F= \D\uo_I$.
According to which is more convenient we can either check if  $A\in
I \Rightarrow f^{-1}[A] \in I$ or if $B \in F \Rightarrow f^{-1}[B]
\in F$.

\medskip

 From now on we consider the clones
$\S_F$ on their own right, not mentioning ideals for a while.

For the two
extreme filters $ \{X\}$ and $\P(X)$ all operations are
continuous. Call a filter {\em proper} iff it is distinct from
the two extreme ones.
In contrast to the ideal case (where the ideal of all finite
sets was an exception) no proper filter yields the full clone.
To see this, denote the filter in question by $F$ and
choose $c\not\in A \in F$.
Then the  constant function with value $c$  is not $F$-
continuous, for $c^{-1}[A] = \emptyset \not\in F$. Hence,
$\S_F^{(1)} \not= \O^{(1)}$.

Different proper filters yield
 clones with different unary parts. To see this, consider $F\not=G$ and choose
(wlog) $B\in F , B \not\in G$ and $c\in   X$ such that $X
\setminus \{c\}\in  G$. Define
$$
    f(x)= \left\{\begin{array}{rl}
    x,& x\in B\\ c,& x\not\in B\end{array}\right. .
$$
Then $f$ belongs
to $\S_F$, for, $A\in F$ implies $f^{-1}[A] \supseteq  A \cap
B\in F$. But $f$ is not $G$-continuous, for $f^{-1}[X\setminus
\{c\}] \subseteq B \not\in G$.

Some clones $\S_F$ come with a handicap that prevents them outright
from being maximal: they are not maximal in their \emph{monoidal
interval}, that is, in the set of clones $\C$ which have the same
unary fragment $\C\uo$ as $\S_F$. To remedy that, we also consider
the biggest clone having $\S\uo_F$ as its unary part. It will be
denoted by $\U_F$ and can
 be described as follows.
$$
    \U_F = \Pol\left(\S_F^{(1)}\right):=\newline
    \bigcup _{n\geq 1} \left\{ f\in
    \O^{(n)}: f(g_1, \dots, g_n) \in \S\uo_F  \mbox{ for all } g_1,
    \dots g_n \in \S_F^{(1)} \right\}.
$$
Then $\U_F^{(1)} = \S_F^{(1)}$ and $\S_F \subseteq \U_F$. In the
dual case the same procedure does not lead anywhere; as  already
mentioned in \ref{ukrit}: $\Pol(\C_I\uo) = \C_I$ (exercise).
\medskip

By now, it is  not well-understood under what conditions $\S_F =
\U_F$. Leaving the proof as an exercise we mention that $\S_F =
\U_F$ holds for countably generated filters and admit that we
do not  know if this also holds in the general regular case. In contrast, we
have the following
\medskip

\begin{prop}\label{continuous_not_equal_for_ultrafilter}
Let $F$ be a non-principal ultrafilter.
\BE
\item  $\U_F = \{ f\in \O: fix(f)\in F\}$,
where $fix(f)$ denotes the fixed-point-set of $f$, i.e., $fix
(f) := \{ x\in X:
f(x,x,\dots, x) = x\}$.
\item $\S_F$ is a proper subset of   $\U_F$.
\EE
\end{prop}

\medskip

The  proof of  (1) is based on the following
result:


\begin{lem}[{\sc Katetov}]
\label{katetov}
A unary  function is $F$-continuous, iff its
fixed-point-set belongs to $F$.
\end{lem}
\begin{proof}[Proof of Lemma~\ref{katetov}.]
The implication $fix(f)\in F \Rightarrow f\in \U_F$ is clear.

Assume that the unary function~$f$ is~$F$-continuous. We want to show
that $fix(f)\in F$.
Let $C:= \{x: \exists k\, f^k(x)=x\}$.
Any unary
function~$f$ defines an undirected graph on~$X$ with edges $(x, f(x))$.  In
each connected component of~$X$ pick a representative --- if possible, in~$C$.
Let~$B$ be the set of those representatives.
Notice that $B \cap f\inv[B] \subseteq fix(f)$.

  For each $x\in X$ let $n(x) $ be minimal such that
$f^n(x)\in B$;  if this is not defined, let $n(x):= \min \{k+j:
\exists \,b\in B \  f^k(x)=f^j(b)\}$. In other words, $n(x)$ is
\begin{enumerate}
\item[--] the length of the unique path from $x$ to an element of $B$, if
$x$ is in a component without fixed points or cycles
\item[--] the smallest $n$ with $f^n(x)\in B$, otherwise.
\end{enumerate}
Let $X_i:=\{x\in X: n(x)\equiv i \pmod 2\}$  for $i=0,1$. It is easy
to see that $f\inv[X_i] \subseteq X_{1-i}\cup B$. Clearly, one of
the $X_i$ is in the ultrafilter $F$ and the other is not. Assume
$X_i\in F$. By $F$-continuity, $f\inv[X_i] \in F$, hence
$X_{1-i}\cup B\in F$. Now, $X_{1-i}\not\in F$ yields $B\in F$. But
then also $f\inv[B]\in F$, so  $F \ni B\cap f\inv[B]\subseteq
fix(f)$.

\end{proof}

\medskip
We now prove Proposition~\ref{continuous_not_equal_for_ultrafilter} for functions
of arbitrary arity:
\begin{proof}[Proof  of Proposition~\ref{continuous_not_equal_for_ultrafilter}]
Let $f\in \U_F$ be given. Then $f(\id,\id,\dots, \id)$ is unary and
$F$-continuous, so its fixed-point-set, which is $fix(f)$ belongs to
$F$.

The other way round. Assuming $fix(f) \in F$ and $g_1, \dots, g_n\in
\S_F^{(1)}$, we must prove $f(g_1, \dots, g_n) \in \S_F^{(1)}$,
i.e.\ $fix \Big( f(g_1, \dots, g_n) \Big) \in F$. But this follows
from the obvious
\[
F \ni fix(f) \cap fix(g_1)\cap \dots \cap fix(g_n) \subseteq
fix \Big( f(g_1, \dots, g_n) \Big).
\]

To exhibit $f\in \U_F \setminus \S_F$, and thus proving (2), we
let $0$ be any point in $X$ and define
$$
    f(x,y) := \left\{\begin{array}{rl} x,& x=y\\
    0,& \mbox{ otherwise}\end{array}\right. .
$$
Then $fix(f) = X \in F$. But for any infinite $B$ we have $f[B^2]
\ni 0$, hence $f[B^2] \subseteq X\setminus \{0\}$ is impossible,
disproving $F$-continuity.
\end{proof}

{\it Remark. } With the description $(1)$ the clones $\U_F$ were
used in \cite{Mar81-CardinalityOfTheSetOfPrecomplete}
by {\sc Marchenkov}, who showed them maximal and distinct for different
ultrafilters. These were the first easy examples of
$2^{2^{\aleph_0}}$ maximal clones.

\medskip

Next we characterize the precomplete  clones of type $\U_F$. Notice
that the following theorem is true without the countability
assumption on $X$ neither do we assume that all cofinite sets are in
$F$. The filter has just to be proper, i.e.\ different from $\{X\}$
and $\P(X)$.

\begin{thm}\label{filtermaxtest}
 If $F$ is a proper filter on $X$, then each of the following
conditions is equivalent to the precompleteness of $\U_F$.
\BE
\item There is no proper filter $G \supsetneq F$ with $\S_F^{(1)}
\subsetneq \S_G^{(1)}$.
\item For each $A\not\in F$ there exists $f\in \S_F^{(1)}$ such
that $f^{-1}[A] = \emptyset$.
\item $\S_F(h) = \O$ for each unary $h\not\in \S_F$.
\EE
\end{thm}

{\it Remarks. }  Condition (1) is formally weaker than
the maximality of
$\U_F$  among filter clones. These conditions may be
equivalent, however.

 From (3) it does not follow that $\S_F$ is maximal, because
there can be a binary function in $\U_F\setminus \S_F$. This
is the case for ultrafilters.
We do not know if $\U_F$ is generated by $\S_F$ and some binary
non-continuous function.

\medskip

\begin{proof}
The necessity of $(1)$ will be established  by constructing a
clone above $\U_F$ from a proper filter $G \supsetneq F$ such that
$\S_G^{(1)} \supsetneq \S_F^{(1)}$.
The first idea is, of course, trying $\U_G$. But this need not
work, so we have to come up with something more tricky.
We consider
\[
\M := \left\{ f \in \O^{(1)}: \bool{f= \tilde f}\in G  \mbox{ for some }
\tilde f \in \S_F^{(1)} \right\},
\]
where $\bool{f=\tilde f}$ denotes the so-called {\em equalizer} $\{
x\in X: f(x) = \tilde f(x)\}$.
A number of easy verifications then yields that
 $\M$ is a monoid and $\S_F^{(1)} \subseteq \M
\subseteq \S_G^{(1)}$.
To see, for example,  that $\M$ is closed under composition, consider  $f,g\in
\M$ witnessed by  $\tilde f, \tilde g \in \S_F^{(1)}$.
Then $\tilde f \circ \tilde g$ witnesses $f\circ g \in \M$, because
\[
\bool{ f\circ g = \tilde f \circ \tilde g} \supseteq \bool{g= \tilde
g} \cap \tilde g^{-1}\left[ \bool{f=\tilde f }\right]
\]
belongs to $G$ (for the preimage
 $\S^{(1)}_F \subseteq \S_G^{(1)}$ is used).
To see that $\M \subseteq \S\uo_G$, let
$f\in \M$ be witnessed by $\tilde f$. Because $\S\uo_F \subseteq
\S\uo_G$, we have $\tilde f \in \S\uo_G$, therefore,
$
f^{-1}[A]   \supseteq \tilde f^{-1}[A] \cap \bool{f=\tilde f}
$
belongs to $G$, whenever $A\in G$.

\medskip

Now we prove $\U_F = \Pol\left(\S\uo_F\right) \subseteq
\Pol\left(\M\right)$. Consider  $h(x_1, \dots, x_n) \in \U_F$ and
$f_1, \dots f_n \in \M$, witnessed by $\tilde f_1, \dots, \tilde f_n
\in \S^{(1)}_F$. Then the function $h(\tilde f_1, \dots,\tilde f_n)
$ witnesses $h(f_1, \dots, f_n) \in \M$. For, $h \in
\Pol\left(\S\uo_F\right) $ implies $h(\tilde f_1, \dots, \tilde
f_n)\in \S\uo_F$, and
\[
\bool{h(f_1, \dots, f_n) = h(\tilde f_1, \dots, \tilde f_n)}
\supseteq \bool{f_1=\tilde f_1} \cap \dots \cap \bool{f_n=\tilde f_n}
\in G.
\]

As $G$ is proper and  $\M \subseteq \S^{(1)}_G \not= \O^{(1)}$ we
cannot have   $\Pol(\M)=\O$. It remains to exhibit a function in
$\Pol(\M)$  which is not in $\U_F$. In fact, we find one in
$\M\setminus \S_F\uo$. Take $A \in G\setminus F$ and $c\in X$ such
that $X \setminus\{c\} \in F$. The former is possible because $G$
strictly includes $F$; the latter because $F$ is proper. Then the
function $p(x):=  \left\{\begin{array}{rl} x,& x\in A\\ c,& x\not\in
A \end{array}\right.$ is in $\M$, because $\bool{p = \id} \supseteq
A \in G$, but not $F$-continuous, because $p^{-1}[ X \setminus
\{c\}] \subseteq A \not\in F$.

\medskip

Next we prove $(1) \Rightarrow (2)$. Actually, we assume that
$A\not\in  F$ is a counterexample to $(2)$, i.e.\ $f^{-1}[A]
\not=\emptyset$ for all $f\in \S\uo_F$, and show that
\[
G:= \{ B \subseteq X: f^{-1}[A] \subseteq B \mbox{ for some }
f\in \S\uo_F\}
\]
is a filter contradicting (1). It is clear that $G$ is
upward-closed.
To see that $G$ is closed under unions, observe that
 $h^{-1}[A] = f_1^{-1}[A] \cup f_2^{-1}[A]$ for the function
$$
    h(x) =  \left\{\begin{array}{rl} f_1(x),& f_1(x)  \in A\\ f_2(x),&
    \mbox{ otherwise} \end{array}\right. .
$$
If $f_{1}, f_2 \in \S\uo_F$,
then so is $h$, for $h^{-1}[B] \supseteq f_1^{-1}[B] \cap
f_2^{-1}[B]\in F$ for any $B\in F$.

 The choice of $A$ just means that $\emptyset
\not\in G$. From $ A\in G$ we conclude that  $G$ is proper and
not equal to $F$.
As $\S\uo_F\subseteq \S\uo_G$ is obvious, it just remains to show
$F\subseteq G$.
Let $B\in F$ be given and let $c$ be some element outside $A$
(if there were none, we had  the impossible $X=A\not\in F$).
Then the function
$$
    q(x) = \left\{\begin{array}{rl} x,& x  \in B\\ c,&
    x\not\in B \end{array}\right.
$$
belongs to $\S\uo_F$ and $q^{-1}[A]
\subseteq B$ proves $B\in G$.

Next we prove $(2) \Rightarrow (3)$. Let $h \not\in \S\uo_F$ be
given. Then $h^{-1}[B] \not\in F$ for some $B\in F$. Using (2)
we take $f\in \S\uo_F$  such that $f^{-1}[h^{-1}[B]] = \emptyset$.
In other words $h \circ f$ does not take any value in $B$.

Our aim is to show $\S_F(h) = \O$. So let any $g(x_1, \dots, x_n)$
be given. Put
$$
    \tilde g(x_1, \dots, x_n, y):= \left\{\begin{array}{rl} g(x_1,
    \dots, x_n),& y\not\in B\\ y,& y \in B. \end{array}\right. .
$$
Then $\tilde g$ is $F$-continuous. For, if
$C\in F$ then $B \cap C\in F$ and $\tilde g[(B\cap C)^{n+1}] = B
\cap C \subseteq C$.

As, obviously, $g(x_1, \dots, x_n) = \tilde g( x_1, \dots, x_n,
h(f(x_1)))$, we have $g \in \S_F(h)$.

It remains to see that (3) is sufficient for the precompleteness
of $\U_F$. For  an arbitrary operation $g(x_1, \dots, x_n)$
outside $\U_F$ we have to show that $\U_F(g) = \O$.

 From $g\not\in \U_F = \Pol(\S\uo_F)$ we get $f_1, \dots, f_n \in
\S\uo_F$ such that $h:= g(f_1, \dots, f_n) \not\in \S\uo_F$.
But then, by (3),
$$
    \O = \S_F(h) \subseteq \S_F(g) \subseteq \U_F(g).
$$
\end{proof}

The theorem is now completely proved. If true, condition (2) is usually
easy to verify. It  yields the precompleteness  of $\U_F$ for, e.g.,
ultrafilters and  countably generated filters.

No new considerations are, however, needed in these examples,
because it is possible to relate the precompleteness of $\U_F$
to that of $\C_I$ for the dual ideal.
This is  at first sight  surprising
 because these clones sit in rather different parts of the lattice.
The  operations of $\U_F$ and $\C_I$ are very
different: For example, $\C_I$ contains all constant operations,
whereas $\U_F$ will never contain any constant operation. In a very
free interpretation, one could say that the unary operations in
$\U_F$ are in a way close to injective, since the preimages of small
sets are small, whereas an operation is more likely to belong to
$\C_I$ the less injective it is.

In the following the countability of $X$ is essential again, and
$I$ must be an ideal `in our sense'.

\begin{cor}\label{precompletnessconnections}
 Let $I$ be  any ideal  and $F$ its dual filter. $\U_F$
is maximal iff $T(A,1,\infty)$ holds for all $A\not\in I$, i.e.\ iff
$\C_I$ is maximal via (possibly infinitely many) unary functions.
\end{cor}

This follows easily from what we have already proved.

If we switch to complements, condition (2) of the last theorem says
that for all $A\not\in I$ there is some $f\in \S\uo_F$ such that
$f^{-1}[A] = X$. The latter amounts to $f: X \to A$ and $f\in
\S\uo_F$ can be read as $f^{-1}[B] \in I$ for all $B\in I$.

In other words condition (2) says the same as the case $n=1$ of
Lemma~\ref{preim}, where  it was
shown to be  equivalent to $ T(A,1,\infty)$ .
\hfill $\Box$

As an immediate consequence of
Corollary~\ref{precompletnessconnections} we have that there exists
an ideal $I$ such that $\C_I$ is maximal while $\U_F$ is not
maximal. Just choose $I$ (using
Example~\ref{maximal_but_not_via_unary} with $n=1$) such that binary
functions are required to check the precompleteness of $\C_I$.

\section{Uncountable base sets}\label{sec:uncountable}

We briefly discuss the possibility of generalizing the results of this paper
to uncountable base sets. As in the countable case, we may assume
that all ideals have full support, that is, they contain all finite
subsets of $X$. For countable $X$, the assumption that an ideal
contains at least one infinite set and does not contain some
infinite set implies that the induced ideal clone is proper, i.e.\
it does not contain all operations on $X$. This is no longer the
case for uncountable $X$. Define for all infinite $\lambda\leq |X|$
an ideal $I_\lambda$ consisting of all sets $S\subseteq X$ with
$|S|<\lambda$. Then we have

\begin{lem}\label{lem:uncount:nontrivialCondition}
    Let $X$ be infinite and let $I$ be a proper ideal with full support.
    Then $\C_I=\O$ iff $I=I_\lambda$ for some infinite $\lambda\leq |X|$.
\end{lem}

    We skip the easy proof.

The preceding lemma immediately makes it  clear that things will be
more complicated for uncountable $X$; in particular, the basic
Proposition~\ref{simpleconnections} does not hold anymore. It can be
replaced by

\begin{prop}\label{prop:uncount:equality}
    Let $I,J$ be ideals such that $\C_I\uo\subseteq\C_J\uo$, and let
    $\lambda\leq |X|$. If $J$ contains a set of size $\lambda$, then
    all sets in $I$ of size at most $\lambda$ are contained in
    $J$.\\
    In particular, if $I$ and $J$ are ideals such that
    $$
        \sup\{|A|:A\in I\}=\sup\{|A|:A\in J\}=:\lambda
    $$
    and
    $$
        \exists A\in I\,(|A|=\lambda)\leftrightarrow \exists A\in
        J\,(|A|=\lambda),
    $$
 then $\C_I=\C_J$ iff $I=J$.
\end{prop}

Again the proof is straightforward, so we skip it, too.

If we demand ideals to contain at least one \emph{large} set, i.e.\
a set of size $|X|$, then Proposition~\ref{simpleconnections} holds,
but still the deeper results of this paper do not generalize,
e.g. the maximality criterion from
\cite{CH01} fails in the following

\begin{ex}\label{prop:uncountable:boundedIdealPrecompleteClone}
    Let  $X$ be  uncountable.  Then there is an  ideal $I$ with
the following properties:
it has full
support and contains a large set; the induced clone
     $\C_I$ is precomplete; but there is some $A\not\in I$ such
that $T(A,n,p)$ fails for all $n$ and $p$.
\end{ex}

We can assume that $X= Y \times \omega$, where $Y$ is
uncountable. This allows us to define the `below'-relation  on $X$ via
$(y_1,n_1) \prec (y_2,n_2) :\iff n_1<n_2$.

Let $I$ denote the ideal of bounded subsets of $X$, i.e.\ $A\in I$
iff there is some $(b,n)$ such that $(a,m) \prec (b,n)$ for all
$(a,m)\in A$.

There are, of course, countable unbounded sets, i.e.\ $\{y_0\}\times
\omega$. But these cannot be mapped onto $X$ by any finitary
function. So the test fails.

But  $\C_I$ is precomplete, anyway. To see this, consider some $f\not\in \C_I$ (wlog
unary). Then $f$ maps a bounded set $S$ to an unbounded set $U$.
 Now let $g\in\O^{(n)}$ be
    arbitrary. Define an operation $h\in\O^{(n+1)}$ as follows:
    $$
        h(x_1,\ldots,x_n,y)=\left\{\begin{array}{rl} g(x_1,\dots,x_n),&
         g(x_1,\ldots,x_n)\prec y\\ y,& \mbox{
        otherwise.}\end{array}\right.
    $$
    Since $h(x_1,\dots,x_n,y)\prec y$ for all $x_1,\dots,x_n,y\in
    X$, we have $f\in\C_I$. Now define another operation $t\in\O^{(n)}$
    such that
    $$
        t(x_1,\dots,x_n)\in  \{ s\in S:
        g(x_1,\ldots,x_n)\prec f(s)\}.
    $$
    Since $t$ has bounded range, we have $t\in\C_I$. But now
    clearly,
    $g(x_1,\dots,x_n)=h(x_1,\dots,x_n,f(t(x_1,\dots,x_n)))\in \C_I(f)$,
    finishing the proof of precompleteness.

\bigskip

In order to generalize our results to uncountable base sets, the
following restriction on ideals proves convenient: Call an ideal
\emph{suitable} iff it contains all small (small $=$ non-large $=$ of cardinality smaller than $X$) sets,
and contains at least one large set but not all sets. When working
with suitable ideals, all results and proofs of this paper
generalize in a straightforward way, except for the construction of
many precomplete clones without the Axiom of Choice
(Theorem~\ref{many_clones_no_choice}). The necessary big almost
disjoint families exist only under
additional assumptions on cardinal arithmetic. To carry out the
generalization of the other results, the definition of $\hat{I}$
must be adjusted as follows: $\hat{I}=\{A\subseteq X: \mbox{
for all large  }
B\subseteq A \mbox{ there exists some large } C\subseteq B
\mbox{ with }C\in I\}$. The corresponding operator $\cdot^\perp$ is the following: For a
family $A\subseteq \P(X)$, $A^\perp :=\{B\subseteq X: \forall C\in A\,
(C\cap B \mbox{ small})\}$.

Observe that whereas the restriction to ideals having full support
is natural and can easily be argued, there is no obvious reason to
consider only suitable ideals, except for them being \ldots
suitable.

\section{Open problems}\label{sec:openproblems}

\begin{prob}
    In Example~\ref{unary_part_smaller_binary_not}, we exhibited two
    incomparable ideal clones with comparable unary fragments. Do
    there exist incomparable ideal clones with comparable $n$-ary
    fragments, where $n>1$?
\end{prob}

\begin{prob}
    Is every ideal clone generated by its binary fragment? (A
    positive answer would yield a negative answer to the previous
    problem for all $n>1$).
\end{prob}

This is even open for many particular ideals.  For example, let
$I_{d=0}$ be the ideal of all subsets $A\subseteq \N$
with upper density 0.   The upper density $\overline d(A)$
is defined as  $\overline d (A)=\overline{\lim}_{n\to\infty} |A\cap
\{0,\ldots, n\}|/(n+1)$.   This ideal is well known and plays
an important role in
analysis and number theory.

\begin{prob} Is the ideal clone $\C_{I_{d=0}}$ generated by its binary
fragment?
\end{prob}

\begin{prob}
Is $\C_{I_{d=0}}$ a precomplete clone?
\end{prob}

\begin{prob}
Is there a prime ideal $P$ on an infinite set $Y$ such for the
ideal  $I = P\times P$ on $X = Y\times Y$ the clone $\C_I$ is
precomplete via infinitely many unary functions?
\end{prob}

\begin{prob}
Find an ideal $I$ such that $\C_I$ is precomplete but for all
$n$ there is $A\not\in I$ such that $T(A,n,1)$ fails.
\end{prob}

\begin{prob}
Which implications hold between
\[
\U_F \subseteq \U_G, \qquad \S_F \subseteq \S_G \quad \mbox{ and
} \quad
\C_I \subseteq \C_J,
\]
where $F$ and $G$ are the dual filters of the ideals $I$ and
$J$, respectively.
\end{prob}

\begin{prob}
Under what conditions $\U_F = \S_F$ holds? Is this true for regular
filters (i.e.\ duals of regular ideals)?
\end{prob}

Finally, we repeat a problem from \cite{CH01}:

\begin{prob}\label{P3} Let $I$ be an ideal such that $\C_I$ is
not precomplete.   Is there an
ideal $J$ such that the clone
$\C_J$ is an upper cover of $\C_I$ in the clone lattice?
\end{prob}


\providecommand{\bysame}{\leavevmode\hbox
to3em{\hrulefill}\thinspace}

\end{document}